\documentclass[global]{svjour}
\usepackage{amssymb}
\usepackage{amsfonts}
\usepackage{amsmath}

\input{tcilatex}

\begin{document}

\title{{\LARGE Existence of equilibrium for an abstract economy with private
information and a countable space of actions}}
\author{Monica Patriche}
\institute{University of Bucharest 
%TCIMACRO{\TeXButton{email}{\email{monica.patriche@yahoo.com}}}%
%BeginExpansion
\email{monica.patriche@yahoo.com}%
%EndExpansion
}
\mail{University of Bucharest, Faculty of Mathematics and Computer Science, 14
Academiei Street, 010014 Bucharest, Romania}
\maketitle

\begin{abstract}
{\small We define the model of an abstract economy with private information
and a countable set of actions. We generalize the H. Yu and Z. Zhang's model
(2007), considering that each agent is characterised by a preference
correspondence instead of having an utility function. \ We establish two
different equilibrium existence results. }
\end{abstract}



\section{\textbf{\ INTRODUCTION}}

We define the model of an abstract economy with private information and a
countable set of actions. The preference correspondences need not to be
represented by utility functions. The equilibrium concept is an extension of
the deterministic equilibrium. We present the H. Yu and Z. Zhang's model in
[18], in which the agents maximize their expected utilities. Our model is a
generalization of H. Yu and Z. Zhang's one.

A purpose in this paper is to prove the existence of equilibrium for an
abstract economy with private information and a countable set of actions.
The assumptions on correspondences refer to upper semicontinuity and
measurability.

The existence of pure strategy equilibrium for a game with finitely many
players, finite action space and diffuse and independent private information
was first proved by Radner and Rosenthal [15]. This result was extended by
Khan and Sun [10] to the case of a finite game with diffuse and independent
private information and with countable compact metric spaces as their action
spaces. These authores have shown in [8] that Radner and Rosenthal's result
can not be extended to a general action space. We quote the papers of M.A.
Khan, K. Rath, Y. Sun [7], [8], [9] and M.A. Khan, Y. Sun [11], [12],
concerning this subject of research. In [19], H. Yu and Z. Zhang showed the
existence of pure strategy equilibrium for games with countable complete
metric spaces and worked with compact-valued correspondences. They relied on
the Bollob\'{a}s and Varopulos's extension [2] of the marriage lemma to
construct a theory of the distribution of a correspondence from an atomless
probability space to a countable complete metric space. They also studied
the case of the game with a continuum of players.

The classical model of Nash [14] was generalized by many authors. Models
were proposed in his pioneering works by Debreu [3] or later by A. Borglin
and H. Keiding [2], Shafer and Sonnenschein [16], Yannelis and Prahbakar
[18]. Yannelis and Prahbakar developed new tehniques of work for showing the
existence of equilibrium. That is the reason for what we defined a new model
that can be integrated in this direction of development of the game theory.
We use the fixed point method of finding the equilibrium, precisely we use
Ky Fan fixed point theorem for upper semicontinuous correspondences.

The paper is organised as follows: In section 2, some notation and
terminological convention are given. In section 3, H. Yu and Z. Zhang's
expected utility model with a finite number of agents and private
information and their main result in [19] are presented. Section 4
introduces our model, that is, {\normalsize an abstract economy with private
information and a countable space of actions. Section 5 contains existence
results for upper semicontinuous correspondences.}

\section{\textbf{PRELIMINARIES AND NOTATION\protect\smallskip \protect%
\medskip }}

Throughout this paper, we shall use the following notations and definitions:

Let $A$ be a subset of a topological space $X$.

\begin{enumerate}
\item $\tciFourier (A)$ denotes the family of all non-empty finite subsets
of $A$.

\item $2^{A}$ denotes the family of all subsets of $A$.

\item cl $A$ denotes the closure of $A$ in $X$.

\item If $A$ is a subset of a vector space, co$A$ denotes the convex hull of 
$A$.

\item If $F$, $G:$ $X\rightarrow 2^{Y}$ are correspondences, then co$G$, cl $%
G$, $G\cap F$ $:$ $X\rightarrow 2^{Y}$ are correspondences defined by $($co$%
G)(x)=$co$G(x)$, $($cl$G)(x)=$cl$G(x)$ and $(G\cap F)(x)=G(x)\cap F(x)$ for
each $x\in X$, respectively.$\medskip $
\end{enumerate}

\textbf{Definition 1. }Each correspondence\textbf{\ }$F:$ $X\rightarrow
2^{Y} $ has two natural inverses:

\begin{enumerate}
\item the upper inverse $F^{u}$ (also called the \textit{strong inverse}) of
a subset $A$ of $Y$ is defined by $F^{u}(A)=\left\{ x\in A:F(x)\subset
A\right\} .$

\item the lower inverse $F^{l}$ (also called the \textit{weak inverse}) of a
subset $A$ of $Y$ is defined by $F^{l}(A)=\left\{ x\in A:F(x)\cap
A\not=\emptyset \right\} .\medskip $
\end{enumerate}

\textbf{Definition 2}. Let $X$, $Y$ be topological spaces and $%
F:X\rightarrow 2^{Y}$ be a correspondence.

1. $F$ is said to be \textit{upper semicontinuous} if for each $x\in X$ and
each open set $V$ in $Y$ with $F(x)\subset V$, there exists an open
neighborhood $U$ of $x$ in $X$ such that $F(y)\subset V$ for each $y\in U$.

2. $F$ is said to be \textit{lower semicontinuous} \textit{(l.s.c)} if for
each $x\in X$ and each open set $V$ in $Y$ with $F(x)\cap V\neq \emptyset $,
there exists an open neighbourhood $U$ of $x$ in $X$ such that $F(y)\cap
V\neq \emptyset $ for each $y\in U$.\medskip

\textbf{Lemma} 1 [20]. \textit{Let }$X$\textit{\ and }$Y$\textit{\ be two
topological spaces and let }$A$\textit{\ be a closed (resp. open) subset of }%
$X.$\textit{\ Suppose }$F_{1}:X\rightarrow 2^{Y}$\textit{, }$%
F_{2}:X\rightarrow 2^{Y}$\textit{\ are lower semicontinuous (resp. upper
semicontinuous) such that }$F_{2}(x)\subset F_{1}(x)$\textit{\ for all }$%
x\in A.$\textit{\ Then the correspondence }$F:X\rightarrow 2^{Y}$\textit{\
defined by}

\begin{center}
$\mathit{F(x)=}\left\{ 
\begin{array}{c}
F_{1}(x)\text{, \ \ \ \ \ \ \ if }x\notin A\text{, } \\ 
F_{2}(x)\text{, \ \ \ \ \ \ \ \ \ if }x\in A%
\end{array}%
\right. $
\end{center}

\textit{is also lower semicontinuous (resp. upper semicontinuous).\medskip }

\textbf{Definition 3 }Let $(T$, $\mathcal{T})$ be a measurable space, $Y$ a
topological space and $F:T\rightarrow 2^{Y}$ a corespondence.

\begin{enumerate}
\item $F$ is \textit{weakly} \textit{measurable }if $F^{l}(A)\in \mathcal{T}$
for each open subset $A$ of $Y;$

\item $F$ is \textit{measurable }if $F^{l}(A)\in \mathcal{T}$ for each
closed subset $A$ of $Y.$
\end{enumerate}

\textit{Remark.} Let $(T$, $\mathcal{T})$ be a measurable space, $Y$ a
countable set and $F:T\rightarrow 2^{Y}$ a corespondence. Then $F$ is
measurable if for each $y\in Y,F^{-1}(y)=\left\{ t\in T:y\in F(t)\right\} $
is $\mathcal{T-}$measurable.\medskip

\textbf{Lemma 2 [1]. }For\textbf{\ }a correspondence $F:T\rightarrow 2^{Y}$
from a measurable space into a metrizable space we have the following:

\begin{enumerate}
\item If $F$ is measurable, then it is also weakly measurable$;$

\item If $F$ is compact valued and weakly measurable, it is measurable.$%
\medskip $
\end{enumerate}

\textbf{Definition 4 [19]. }Let\textbf{\ }$Y$ be a countable complete metric
space, $(T$, $\mathcal{T}$, $\lambda )$ an atomless probability space and $%
F:T\rightarrow 2^{Y}$ a measurable corespondence. The function $%
f:T\rightarrow Y$ is said to be a selection of $F$ if $f(t)\in F(t)$ for $%
\lambda -$almost $t\in T.$ Denote $\mathcal{D}_{F}=\left\{ \lambda f^{-1}:f%
\text{ is a measurable selection of }F\right\} .\medskip $

\textbf{Lemma 3 [19]. }Let\textbf{\ }$Y$ be a countable complete metric
space, $(T$, $\mathcal{T}$, $\lambda )$ an atomless probability space and $%
F:T\rightarrow 2^{Y}$ a measurable corespondence. Then $\mathcal{D}_{F}$ is
nonempty and convex in the space $\mathcal{M}(Y)$ - the space of probability
measure on $Y$, equipped with the topology of weak convergence.\medskip

\textbf{Lemma 4 [19]}. Let\textbf{\ }$Y$ be a countable complete metric
space, and $(T$, $\mathcal{T}$, $\lambda )$ be an atomless probability space
and $F:T\rightarrow 2^{Y}$ be a measurable corespondence. If $F$ is compact
valued, then $\mathcal{D}_{F}$ is compact in $\mathcal{M}(Y).\medskip $

\textbf{Lemma 5 [19].} Let $X$ be a metric space, $(T$, $\mathcal{T}$, $%
\lambda )$ be an atomless probability space, $Y$ be a countable complete
metric space and $F:T\times X\rightarrow 2^{Y}$ a correspondence. Assume
that for any fixed $x$ in $X$, $F(\cdot ,x)$ (also denoted by $F_{x})$ is a
compact-valued measurable correspondence$,$ and for each fixed $t\in T,$ $%
F(t,\cdot )$ is upper semicontinuous on $X$. Also, assume that there exists
a compact valued corespondence $H:T\times X\rightarrow 2^{Y}$ such that $%
F(t,x)\subset H(t)$ for all $t$ and $x$. Then $\mathcal{D}_{F_{x}}$ is upper
semicontinuous on $X.$

\textbf{Theorem} 1 \textbf{(Kuratowski-Ryll-Nardzewski Selection \ Theorem)
[1]. }A weakly measurable correspondence with nonempty closed values from a
measurable space into a Polish space admits a measurable selector.\medskip

\section{A Nash Equilibrium Existence Theorem}

We present Yu and Zhang's model of a finite game with private information.
In this model it is assigned to each agent a private information related to
his action and payoff described by the random mappings $\tau _{i}$ and $\chi
_{i},$ mappings defined on $(\Omega ,\mathcal{F})=\underset{i\in I}{(\prod }%
(Z_{i},X_{i}),\underset{i\in I}{\prod }(\mathcal{Z}_{i},\mathcal{X}_{i})),$
where $(X_{i}$, $\mathcal{X}_{i})$ and $(Z_{i},\mathcal{Z}_{i})$ are
measurable spaces. For a point $\omega =(z_{1},x_{1},...,z_{n},x_{n})\in
\Omega ,$ $\tau _{i}$ and $\chi _{i}$ are the coordinate projections

$\tau _{i}(\omega )=z_{i}$, $\chi _{i}(\omega )=x_{i}.$ Each player $i$ in $%
I $ first observes the realization, say $z_{i}\in Z_{i},$ of the random
element $\tau _{i}(\omega ),$ then chooses his own action from a nonempty
compact subset $D_{i}(z_{i})$ of a countable complete metric space $A_{i},$
with $D_{i}(\cdot )$ measurable. The payoff of each player $i$ is given by
the the utility function $u_{i}:A\times X_{i}\rightarrow \mathbb{R},$ where $%
A=\underset{j\in I}{\prod }A_{j}$ is the set of of all combinations of all
players' moves. Let $\mu $ be a probability measure on $\Omega .$ It is
assumed the following uniform integrability condition (UI):

(UI) For every $i\in I,$ there is a real-valued integrable function $h_{i}:$ 
$\Omega \rightarrow \mathbb{R}$ such that $\mu -$almost all $\omega \in
\Omega ,$ $\mid u_{i}(a,\chi _{i}(\omega ))\mid \leq h_{i}(\omega )$ holds
for all $a\in A.$

\textbf{Definition 5 [19]}. A \textit{finite game with private information}
is a family $\Gamma =(I,((Z_{i},\mathcal{Z}_{i}),(X_{i},\mathcal{X}%
_{i}),(A_{i},D_{i}),u_{i})_{i\in I},\mu )$.\medskip

For each $i\in I,$ let meas($Z_{i},D_{i})$ be the set of measurable mappings 
$f$ from $Z_{i}$ to $A_{i}$ such that $f(z_{i})\in D_{i}(z_{i})$ for each $%
z_{i}\in Z_{i}.$ An element $g_{i}$ of meas($Z_{i},A_{i})$ is called \textit{%
a pure strategy} for player $i.$ \textit{A pure strategy profile} $g$ is an
n-vector function $(g_{1},g_{2},...,g_{n})$ that specifies a pure strategy
for each player.\medskip

\textbf{Definition 6 [19]}. For a pure strategy profile $%
g=(g_{1},g_{2},...,g_{n}),$ the expected payoff for player $i$ is

$U_{i}(g)=\tint\limits_{\omega \in \Omega }u_{i}(g_{1}(\tau _{1}(\omega
)),...,g_{n}(\tau _{n}(\omega )),\chi _{i}(\omega ))\mu d(\omega ).\medskip $

\textbf{Definition 7 [19].} \textit{A Nash equilibrium in pure strategies}
is defined as a pure strategy profile $(g_{1}^{\ast },g_{2}^{\ast
},...,g_{n}^{\ast })$ such that for each player $i\in I$

$U_{i}(g^{\ast })\geq U_{i}(g_{i},g_{-i}^{\ast })$ for all $g_{i}\in $Meas$%
(Z_{i},D_{i}).\medskip $

The following theorem is the main result of Zhang in [19].

\textbf{Theorem 2\ }\textit{Suppose that for every player }$i,$\textit{\ the
compact valued }$D_{i}$\textit{\ corespondence is measurable, and}

\textit{a) the distribution }$\mu \tau _{i}^{-1}$\textit{\ of }$\tau _{i}$%
\textit{\ is an atomless measure;}

\textit{b) the random elements }$\left\{ \tau _{j}:j\not=i\right\} $\textit{%
\ together with the random element }$\xi _{i}\equiv (\tau _{i},\chi _{i})$%
\textit{\ form a mutually independent set;}

\textit{c) for any fixed }$x_{i}\in X_{i},$\textit{\ }$u_{i}(\cdot ,x_{i})$%
\textit{\ is a continuous function on }$A;$\textit{\ for any fixed }$a\in A,$%
\textit{\ }$u_{i}(a,\cdot )$\textit{\ is a measurable function on }$(X_{i},%
\mathcal{X}_{i});$

\textit{d) the uniform integrability condition (UI) holds.}

\textit{Then the game }$\Gamma $\textit{\ has a Nash equilibrium in pure
strategies.\medskip }

\section{The Model of an abstract economy with private information}

In this section we define a model of abstract economy with private
information and a countable set of actions. We also prove the existence of
equilibrium of abstract economies.

Let $I$ be a nonempty and finite set (the set of agents). For each $i\in I$,
the space of actions, $A_{i}$ is a countable complete metric space and $%
(Z_{i},\mathcal{Z}_{i})$ is measurable space. Let $(\Omega ,\mathcal{F})$ be
the product measurable space$,\underset{i\in I}{(\prod }Z_{i},\underset{i\in
I}{\prod }\mathcal{Z}_{i})$, and $\mu $ a probability measure on $(\Omega ,%
\mathcal{F}).$ For a point $\omega =(z_{1},...,z_{n})\in \Omega ,$ define
the coordinate projections

$\tau _{i}(\omega )=z_{i}.$

The random mapping $\tau _{i}(\omega )$ is interpreted as player i's private
information related to his action.

For each $i\in I,$ we also denote by meas($Z_{i},A_{i})$ the set of
measurable mappings $f$ from $Z_{i}$ to $A_{i}.$ An element $g_{i}$ of meas($%
Z_{i},A_{i})$ is called a \textit{pure strategy} for player $i.$ A \textit{%
pure strategy profile} $g$ is an n-vector function $(g_{1},g_{2},...,g_{n})$
that specifies a pure strategy for each player.

We suppose that there exists a correspondence $D_{i}:Z_{i}\rightarrow
2^{A_{i}}$ such that each agent $i$ can choose an action from $%
D_{i}(z_{i})\subset A_{i}$ for each $z_{i}\in Z_{i}.$

Let $D_{\mathcal{D}_{i}}$\ be the set $\left\{ (\mu \tau
_{i}^{-1})g_{i}^{-1}:g_{i}\text{ is a measurable selection of }D_{i}\right\}
.$

For each $i\in I,$ let the constraint correspondence be $\alpha
_{i}:Z_{i}\times \tprod\limits_{i\in I}D_{D_{i}}\rightarrow 2^{A_{i}}$, such
that $\alpha _{i}(z_{i},(\mu \tau _{1}^{-1})g_{1}^{-1},(\mu \tau
_{2}^{-1})g_{2}^{-1},...,(\mu \tau _{n}^{-1})g_{n}^{-1})\subset A_{i}$ and
the preference correspondence is $P_{i}:Z_{i}\times \tprod\limits_{i\in
I}D_{D_{i}}\rightarrow 2^{A_{i}}$, such that $P_{i}(z_{i},(\mu \tau
_{1}^{-1})g_{1}^{-1},$\newline
$(\mu \tau _{2}^{-1})g_{2}^{-1},...,(\mu \tau _{n}^{-1})g_{n}^{-1})\subset
A_{i}.$\textit{\ }$\medskip $

\textbf{Definition 8.}\textit{\ }An \textit{abstract economy} \textit{(or a
generalized game)} \textit{with private information} \textit{and a countable
space of actions} is defined as $\Gamma =(I,((Z_{i},\mathcal{Z}_{i}),\newline
(A_{i},\alpha _{i},P_{i}))_{i\in I},\mu )$.\medskip 

\textbf{Definition 9}.\textit{\ }An \textit{equilibrium} for $\Gamma $ is
defined as a strategy profile $(g_{1}^{\ast },g_{2}^{\ast },...,g_{n}^{\ast
})\in \tprod\limits_{i\in I}$Meas$(Z_{i},D_{i})$ such that for each $i\in I:$

1) $g_{i}^{\ast }(z_{i})\in \alpha _{i}(z_{i},(\mu \tau
_{1}^{-1})(g_{1}^{\ast })^{-1},(\mu \tau _{2}^{-1})(g_{2}^{\ast
})^{-1},...,(\mu \tau _{n}^{-1})(g_{n}^{\ast })^{-1})$

for each $z_{i}\in Z_{i};$

2) $\alpha _{i}(z_{i},(\mu \tau _{1}^{-1})(g_{1}^{\ast })^{-1},(\mu \tau
_{2}^{-1})(g_{2}^{\ast })^{-1},...,(\mu \tau _{n}^{-1})(g_{n}^{\ast
})^{-1})\cap $

$P_{i}(z_{i},(\mu \tau _{1}^{-1})(g_{1}^{\ast })^{-1},(\mu \tau
_{2}^{-1})(g_{2}^{\ast })^{-1},...,(\mu \tau _{n}^{-1})(g_{n}^{\ast
})^{-1})=\phi $ for each $z_{i}\in Z_{i}.$

.

\section{Existence of equilibrium for abstract economies\textit{\ }with
private information}

We state some new equilibrium existence theorems for abstract economies.

Theorem 3 is an existence theorem of equilibrium for an abstract economy
with upper semicontinuous correspondences $\alpha _{i}$ and $P_{i}.$\medskip

\textbf{Theorem 3}. \textit{Let\ }$\Gamma =(I,((Z_{i},\mathcal{Z}%
_{i}),(A_{i},\alpha _{i},P_{i}))_{i\in I},\mu )$\textit{\ be an abstract
economy with private information and a countable space of action, where }$I$%
\textit{\ is a finite index set such that for each }$i\in I:$

\textit{a) }$A_{i}$\textit{\ is a countable complete metric space\ and }$%
(Z_{i},\mathcal{Z}_{i})$\textit{\ is a measurable space; }$(\Omega ,F)$%
\textit{\ is the product measurable space }$\underset{i\in I}{(\prod }(Z_{i},%
\mathcal{Z}_{i}))$\textit{\ and }$\mu $\textit{\ an atomless probability
measure on }$(\Omega ,F);$

\textit{b) the correspondence }$D_{i}:Z_{i}\rightarrow 2^{A_{i}}$\textit{\
is measurable with compact values; }

\textit{c) the correspondence }$\alpha _{i}:Z_{i}\times \tprod\limits_{i\in
I}\mathcal{D}_{D_{i}}\rightarrow 2^{A_{i}}$\textit{\ is measurable with
respect to }$z_{i}$\textit{\ and, for all }$z_{i}\in Z_{i},$\textit{\ }$%
\alpha _{i}(z_{i},\cdot ,\cdot ,...,\cdot )$\textit{\ is upper
semicontinuous with nonempty, compact values;}

\textit{d) the correspondence }$P_{i}:Z_{i}\times \tprod\limits_{i\in I}%
\mathcal{D}_{D_{i}}\rightarrow 2^{A_{i}}$\textit{\ is measurable with
respect to }$z_{i}$\textit{\ and, for all }$z_{i}\in Z_{i},$\textit{\ }$%
P_{i}(z_{i},\cdot ,\cdot ,...,\cdot )$\textit{\ is upper semicontinuous with
nonempty, compact values;}

\textit{e) for each }$z_{i}\in Z_{i}$\textit{\ and each }$%
(g_{1},g_{2},...,g_{n})\in \tprod\limits_{i\in I}$\textit{Meas}$%
(Z_{i},A_{i}),$

\textit{\ }$g_{i}(z_{i})\not\in P_{i}(z_{i},(\mu \tau
_{1}^{-1})g_{1}^{-1},(\mu \tau _{2}^{-1})g_{2}^{-1},...,(\mu \tau
_{n}^{-1})g_{n}^{-1});$

f) the set $U_{i}:=$\newline
$\left\{ (z_{i},\lambda _{1},\lambda _{2},...,\lambda _{n})\in Z_{i}\times
\tprod\limits_{i\in I}\mathcal{D}_{D_{i}}:(\alpha _{i}\cap
P_{i})(z_{i},\lambda _{1},\lambda _{2},...,\lambda _{n})=\emptyset \right\} $
is open for each $z_{i}\in Z_{i}$.

\textit{Then there exists }$(g_{1}^{\ast },g_{2}^{\ast },...,g_{n}^{\ast
})\in \tprod\limits_{i\in I}$\textit{Meas}$(Z_{i},A_{i})$\textit{\ an
equilibrium for }$\Gamma .$\textit{\bigskip }

\textit{Proof.} By Lemma 3, $D_{D_{i}}$ is nonempty and convex. By Lemma 4, $%
D_{D_{i}}$ is compact$.$ For each $i\in I$ the set

$U_{i}:=\left\{ (z_{i},\lambda _{1},\lambda _{2},...,\lambda _{n})\in
Z_{i}\times \tprod\limits_{i\in I}\mathcal{D}_{D_{i}}:(\alpha _{i}\cap
P_{i})(z_{i},\lambda _{1},\lambda _{2},...,\lambda _{n})=\emptyset \right\} $
is open

and we define $F_{i}:Z_{i}\times \tprod\limits_{i\in I}\mathcal{D}%
_{D_{i}}\rightarrow 2^{A_{i}}$ by

$F_{i}(z_{i},\lambda _{1},\lambda _{2},...,\lambda _{n})=\left\{ 
\begin{array}{c}
(\alpha _{i}\cap P_{i})(z_{i},\lambda _{1},\lambda _{2},...,\lambda _{n})%
\text{ if }(z_{i},\lambda _{1},\lambda _{2},...,\lambda _{n})\not\in U_{i},
\\ 
\alpha _{i}(z_{i},\lambda _{1},\lambda _{2},...,\lambda _{n})\text{ if }%
(z_{i},\lambda _{1},\lambda _{2},...,\lambda _{n})\in U_{i}.%
\end{array}%
\right. $

Then the correspondence $F_{i}$ has nonempty, compact values and is
measurable with respect to $z_{i}$ and upper semicontinuous with respect to%
\textit{\ }$(\lambda _{1},\lambda _{2},...,\lambda _{n})\in
\tprod\limits_{i\in I}\mathcal{D}_{D_{i}}.$

We denote $\mathcal{D}_{F_{i}}(\lambda _{1},\lambda _{2},...,\lambda _{n})=$

=$\left\{ (\mu \tau _{i}^{-1})g_{i}^{-1}:g_{i}\text{ is a measurable
selection of }F_{i}(\cdot ,\lambda _{1},\lambda _{2},...,\lambda
_{n})\right\} .$ Then:

i) $\mathcal{D}_{F_{i}}(\lambda _{1},\lambda _{2},...,\lambda _{n})$ is
nonempty because there exists a measurable selection from the correspondence 
$F_{i}$ by Kuratowski-Ryll-Nardewski Selection Theorem.

ii) $\ \mathcal{D}_{F_{i}}(\lambda _{1},\lambda _{2},...,\lambda _{n})$ is
convex and compact by Lemma 3 and Lemma 4.

We define $\Phi :$ $\tprod\limits_{i\in I}\mathcal{D}_{D_{i}}\rightarrow
2^{\tprod\limits_{i\in I}\mathcal{D}_{D_{i}}},$ $\Phi (\lambda _{1},\lambda
_{2},...,\lambda _{n})=\tprod\limits_{i\in I}\mathcal{D}_{F_{i}}(\lambda
_{1},\lambda _{2},...$\newline
$,\lambda _{n}).$

The set $\tprod\limits_{i\in I}\mathcal{D}_{D_{i}}$ is nonempty, compact and
convex. By Lemma 5 the correspondence $\mathcal{D}_{F_{i}}$ is upper
semicontinuous. Then the correspondence $\Phi $ is upper semicontinuous and
has nonempty compact and convex values. By Ky Fan fixed point Theorem, we
know that there exists a fixed point $(\lambda _{1}^{\ast },\lambda
_{2}^{\ast },...,\lambda _{n}^{\ast })\in \Phi (\lambda _{1}^{\ast },\lambda
_{2}^{\ast },...,\lambda _{n}^{\ast }).$ In particular, for each player $i,$ 
$\lambda _{i}^{\ast }\in \mathcal{D}_{F_{i}}(\lambda _{1}^{\ast },\lambda
_{2}^{\ast },...,\lambda _{n}^{\ast }).$ Therefore, for each player $i,$
there exists $g_{i}^{\ast }\in $\textit{Meas}$(Z_{i},A_{i})$ such that $%
g_{i}^{\ast }$ is a selection of $F_{i}(\cdot ,\lambda _{1}^{\ast },\lambda
_{2}^{\ast },...,\lambda _{n}^{\ast })$ and $(\mu \tau
_{i}^{-1})(g_{i}^{\ast })^{-1}=\lambda _{i}^{\ast }.$

We prove that $(g_{1}^{\ast },g_{2}^{\ast },...,g_{n}^{\ast })$ is an
equilibrium for $\Gamma .$ For each $i\in I,$ because $g_{i}^{\ast }$ is a
selection of $F_{i}(\cdot ,\lambda _{1}^{\ast },\lambda _{2}^{\ast
},...,\lambda _{n}^{\ast })$, it follows that $g_{i}^{\ast }(z_{i})\in
(\alpha _{i}\cap P_{i})(z_{i},\lambda _{1}^{\ast },\lambda _{2}^{\ast
},...,\lambda _{n}^{\ast })$ if $(z_{i},\lambda _{1}^{\ast },\lambda
_{2}^{\ast },...,\lambda _{n}^{\ast })\not\in U_{i}$ or $g_{i}^{\ast
}(z_{i})\in \alpha _{i}^{^{\prime }}(z_{i},\lambda _{1}^{\ast },\lambda
_{2}^{\ast },...,\lambda _{n}^{\ast })$ if $(z_{i},\lambda _{1}^{\ast
},\lambda _{2}^{\ast },...,\lambda _{n}^{\ast })\in U_{i}.$

By the assumption d) it follows that $g_{i}^{\ast }(z_{i})\not\in
P_{i}(z_{i},\lambda _{1}^{\ast },\lambda _{2}^{\ast },...,\lambda _{n}^{\ast
})$ for each $z_{i}\in Z_{i}.$ Then $g_{i}^{\ast }(z_{i})\in \alpha
_{i}(z_{i},\lambda _{1}^{\ast },\lambda _{2}^{\ast },...,\lambda _{n}^{\ast
})$ and $(z_{i},\lambda _{1}^{\ast },\lambda _{2}^{\ast },...,\lambda
_{n}^{\ast })\in U_{i}.$ This is equivalent with the fact that $g_{i}^{\ast
}(z_{i})\in \alpha _{i}(z_{i},(\mu \tau _{1}^{-1})(g_{1}^{\ast })^{-1},(\mu
\tau _{2}^{-1})(g_{2}^{\ast })^{-1},...,$\newline
$(\mu \tau _{n}^{-1})(g_{n}^{\ast })^{-1})$ and $(\alpha _{i}\cap
P_{i})(z_{i},(\mu \tau _{1}^{-1})(g_{1}^{\ast })^{-1},(\mu \tau
_{2}^{-1})(g_{2}^{\ast })^{-1},...,(\mu \tau _{n}^{-1})(g_{n}^{\ast
})^{-1})=\emptyset $ for each $z_{i}\in Z_{i}.$ Consequently, $(g_{1}^{\ast
},g_{2}^{\ast },...,g_{n}^{\ast })$ is an equilibrium for $\Gamma .\medskip $

\textbf{Theorem 4}. \textit{Let\ }$\Gamma =(I,((Z_{i},\mathcal{Z}%
_{i}),(A_{i},\alpha _{i},P_{i}))_{i\in I},\mu )$\textit{.\ be an abstract
economy with private information and a countable space of action, where }$I$%
\textit{\ is a finite index set such that for each }$i\in I,$

\textit{a) }$A_{i}$\textit{\ is a countable complete metric space and }$%
(Z_{i}$\textit{, }$\mathcal{Z}_{i})$\textit{\ is a measurable space; }$%
(\Omega ,\mathcal{F})$\textit{\ is the product measurable space }$\underset{%
i\in I}{\prod }(Z_{i},\mathcal{Z}_{i})$\textit{\ and }$\mu $\textit{\ an
atomless probability measure on }$(\Omega ,\mathcal{F});$

\textit{b) the correspondence }$D_{i}:Z_{i}\rightarrow 2^{A_{i}}$\textit{\
is measurable with compact values; let }$D_{D_{i}}$\textit{\ be the set}

\textit{\ }$\left\{ (\mu \tau _{i}^{-1})g_{i}^{-1}:g_{i}\text{ \textit{is a
measurable selection of} }D_{i}\right\} ;$

\textit{c) the correspondence }$\alpha _{i}:Z_{i}\times \tprod\limits_{i\in
I}\mathcal{D}_{D_{i}}\rightarrow 2^{A_{i}}$\textit{\ is measurable with
respect to }$z_{i}$\textit{\ and, for all }$z_{i}\in Z_{i},$\textit{\ }$%
\alpha _{i}(z_{i},\cdot ,\cdot ,...,\cdot )$\textit{\ is upper
semicontinuous with nonempty, compact values;}

\textit{d) there exists a selector }$G_{i}:Z_{i}\times \tprod\limits_{i\in I}%
\mathcal{D}_{D_{i}}\rightarrow 2^{A_{i}}$\textit{\ for }$(\alpha _{i}\cap
P_{i}):Z_{i}\times \tprod\limits_{i\in I}\mathcal{D}_{D_{i}}\rightarrow
2^{A_{i}}$ \textit{such that }$G_{i}(z_{i},(\mu \tau
_{1}^{-1})g_{1}^{-1},(\mu \tau _{2}^{-1})g_{2}^{-1},...,(\mu \tau
_{n}^{-1})g_{n}^{-1})$\textit{\ is measurable with respect to }$z_{i}$%
\textit{\ and, for all }$z_{i}\in Z_{i},$\textit{\ }$G_{i}(z_{i},\cdot
,\cdot ,...,\cdot )$\textit{\ is upper semicontinuous with nonempty, compact
values;}

\textit{e) for each }$z_{i}\in Z_{i}$\textit{\ and each }$%
(g_{1},g_{2},...,g_{n})\in \tprod\limits_{i\in I}$\textit{Meas}$%
(Z_{i},A_{i}),$

\textit{\ }$g_{i}(z_{i})\not\in G_{i}(z_{i},(\mu \tau
_{1}^{-1})g_{1}^{-1},(\mu \tau _{2}^{-1})g_{2}^{-1},...,(\mu \tau
_{n}^{-1})g_{n}^{-1});$

f) the set\newline
$U_{i}:=\left\{ (z_{i},\lambda _{1},\lambda _{2},...,\lambda _{n})\in
Z_{i}\times \tprod\limits_{i\in I}\mathcal{D}_{D_{i}}:(\alpha _{i}\cap
P_{i})(z_{i},\lambda _{1},\lambda _{2},...,\lambda _{n})=\emptyset \right\} $
is open for each $z_{i}\in Z_{i}$.

\textit{Then there exists }$(g_{1}^{\ast },g_{2}^{\ast },...,g_{n}^{\ast
})\in \tprod\limits_{i\in I}$\textit{Meas}$(Z_{i},A_{i})$\textit{\ an
equilibrium for }$\Gamma .$\textit{\bigskip }

\textit{Proof.} By Lemma 3, $D_{D_{i}}$ is nonempty and convex. By Lemma 4, $%
\mathcal{D}_{D_{i}}$ is compact$.$ For each $i\in I$ the set

$U_{i}:=\left\{ (z_{i},\lambda _{1},\lambda _{2},...,\lambda _{n})\in
Z_{i}\times \tprod\limits_{i\in I}\mathcal{D}_{D_{i}}:(\alpha _{i}\cap
P_{i})(z_{i},\lambda _{1},\lambda _{2},...,\lambda _{n})=\emptyset \right\} $
is open

and we define $F_{i}:Z_{i}\times \tprod\limits_{i\in I}\mathcal{D}%
_{D_{i}}\rightarrow 2^{A_{i}}$ by

$F_{i}(z_{i},\lambda _{1},\lambda _{2},...,\lambda _{n})=\left\{ 
\begin{array}{c}
G_{i}(z_{i},\lambda _{1},\lambda _{2},...,\lambda _{n})\text{ if }%
(z_{i},\lambda _{1},\lambda _{2},...,\lambda _{n})\not\in U_{i}, \\ 
\alpha _{i}(z_{i},\lambda _{1},\lambda _{2},...,\lambda _{n})\text{ if }%
(z_{i},\lambda _{1},\lambda _{2},...,\lambda _{n})\in U_{i}.%
\end{array}%
\right. $

Then the correspondence $F_{i}$ has nonempty, compact values and is
measurable with respect to $z_{i}$ and upper semicontinuous with respect to%
\textit{\ }$(\lambda _{1},\lambda _{2},...,\lambda _{n})\in
\tprod\limits_{i\in I}\mathcal{D}_{D_{i}}.$

We denote $\mathcal{D}_{F_{i}}(\lambda _{1},\lambda _{2},...,\lambda _{n})=$%
\newline

$=\left\{ (\mu \tau _{i}^{-1})g_{i}^{-1}:g_{i}\text{ is a measurable
selection of }F_{i}(\cdot ,\lambda _{1},\lambda _{2},...,\lambda
_{n})\right\} .$ Then:

i) $\mathcal{D}_{F_{i}}(\lambda _{1},\lambda _{2},...,\lambda _{n})$ is
nonempty because there exists a measurable selection from the correspondence 
$F_{i}$ by Kuratowski-Ryll-Nardewski Selection Theorem.

ii) $\ \mathcal{D}_{F_{i}}(\lambda _{1},\lambda _{2},...,\lambda _{n})$ is
convex and compact by Lemma 3 and Lemma 4.

We define $\Phi :$ $\tprod\limits_{i\in I}\mathcal{D}_{D_{i}}\rightarrow
2^{\tprod\limits_{i\in I}\mathcal{D}_{D_{i}}},$ $\Phi (\lambda _{1},\lambda
_{2},...,\lambda _{n})=\tprod\limits_{i\in I}\mathcal{D}_{F_{i}}(\lambda
_{1},\lambda _{2},...,\lambda _{n}).$

The set $\tprod\limits_{i\in I}\mathcal{D}_{D_{i}}$ is nonempty, compact
convex. By Lemma 5 the correspondence $\mathcal{D}_{F_{i}}$ is upper
semicontinuous. Then the correspondence $\Phi $ is upper semicontinuous and
has nonempty, compact and convex values. By Ky Fan fixed point theorem, we
know that there exists a fixed point $(\lambda _{1}^{\ast },\lambda
_{2}^{\ast },...,\lambda _{n}^{\ast })\in \Phi (\lambda _{1}^{\ast },\lambda
_{2}^{\ast },...,\lambda _{n}^{\ast }).$ In particular, for each player $i,$ 
$\lambda _{i}^{\ast }\in \mathcal{D}_{F_{i}}(\lambda _{1}^{\ast },\lambda
_{2}^{\ast },...,\lambda _{n}^{\ast }).$ Therefore, for each player $i,$
there exists $g_{i}^{\ast }\in $\textit{Meas}$(Z_{i},A_{i})$ such that $%
g_{i}^{\ast }$ is a selection of $F_{i}(\cdot ,\lambda _{1}^{\ast },\lambda
_{2}^{\ast },...,\lambda _{n}^{\ast })$ and $(\mu \tau
_{i}^{-1})(g_{i}^{\ast })^{-1}=\lambda _{i}^{\ast }.$

We prove that $(g_{1}^{\ast },g_{2}^{\ast },...,g_{n}^{\ast })$ is an
equilibrium for $\Gamma .$ For each $i\in I,$ because $g_{i}^{\ast }$ is a
selection of $F_{i}(\cdot ,\lambda _{1}^{\ast },\lambda _{2}^{\ast
},...,\lambda _{n}^{\ast })$, it follows that $g_{i}^{\ast }(z_{i})\in
(\alpha _{i}\cap P_{i})(z_{i},\lambda _{1}^{\ast },\lambda _{2}^{\ast
},...,\lambda _{n}^{\ast })$ if $(z_{i},\lambda _{1}^{\ast },\lambda
_{2}^{\ast },...,\lambda _{n}^{\ast })\not\in U_{i}$ or $g_{i}^{\ast
}(z_{i})\in \alpha _{i}(z_{i},\lambda _{1}^{\ast },\lambda _{2}^{\ast
},...,\lambda _{n}^{\ast })$ if $(z_{i},\lambda _{1}^{\ast },\lambda
_{2}^{\ast },...,\lambda _{n}^{\ast })\in U_{i}.$

By the hypothesis d), it follows that $g_{i}^{\ast }(z_{i})\not\in
G_{i}(z_{i},\lambda _{1}^{\ast },\lambda _{2}^{\ast },...,\lambda _{n}^{\ast
})$ for each $z_{i}\in Z_{i}.$ Then $g_{i}^{\ast }(z_{i})\in \alpha
_{i}(z_{i},\lambda _{1}^{\ast },\lambda _{2}^{\ast },...,\lambda _{n}^{\ast
})$ and $(z_{i},\lambda _{1}^{\ast },\lambda _{2}^{\ast },...,\lambda
_{n}^{\ast })\in U_{i}.$ This is equivalent with the fact that $g_{i}^{\ast
}(z_{i})\in \alpha _{i}(z_{i},(\mu \tau _{1}^{-1})(g_{1}^{\ast })^{-1},(\mu
\tau _{2}^{-1})(g_{2}^{\ast })^{-1},...\newline
,(\mu \tau _{n}^{-1})(g_{n}^{\ast })^{-1})$ and $(\alpha _{i}\cap
P_{i})(z_{i},(\mu \tau _{1}^{-1})(g_{1}^{\ast })^{-1},(\mu \tau
_{2}^{-1})(g_{2}^{\ast })^{-1},...,(\mu \tau _{n}^{-1})(g_{n}^{\ast
})^{-1})=\emptyset $ for each $z_{i}\in Z_{i}.$ Consequently, $(g_{1}^{\ast
},g_{2}^{\ast },...,g_{n}^{\ast })$ is an equilibrium for $\Gamma .\medskip $

\end{document}